\documentclass{article}

\usepackage[dvips]{graphicx}
\usepackage{amsmath,amsfonts,amssymb}

\newtheorem{Definition}{Definition}[section]

\newtheorem{Theorem}[Definition]{Theorem}
\newtheorem{Lemma}[Definition]{Lemma}
\newtheorem{Corollary}[Definition]{Corollary}
\newtheorem{Remark}[Definition]{Remark}

\newcommand{\ed}{{\rm End\,}}
\newcommand{\rad}{{\rm rad\,}}
\newcommand{\soc}{{\rm soc\,}}
\newcommand{\pd}{{\rm proj.dim\, }}
\newcommand{\rp}{{\rm rep.dim\, }}
\newcommand{\gd}{{\rm gl.dim\, }}
\newcommand{\add}{{\rm add \, }}
\newcommand{\Hm}{{\rm Hom \, }}
\newcommand{\repdim}{{\rm repdim \, }}
\newcommand{\End}{{\rm End \, }}

\title{{\sc On the representation dimension of rank 2 group algebras and 
related algebras}}
\date{}
\author{
{\sc Thorsten Holm}\\
University of Leeds\\
\and {\sc Wei Hu}\\
University of Leeds/\\ Beijing Normal University
}

\begin{document}

\maketitle

\begin{abstract}
The representation dimension was defined by M.\,Auslander in 1970
and is, due to spectacular recent progress, one of the most 
interesting homological invariants in representation theory. 
The precise value is not known in general, and is very hard to
compute even for small examples. 
For group algebras, it is known in the case of
cyclic Sylow subgroups, due to Auslander's fundamental work. 
For some group algebras (in characteristic 2) of rank at least 3
the precise value of the representation dimension follows from
recent work of R.\,Rouquier. There is a gap for group
algebras of rank 2; here the deep geometric methods
do not work. In this paper we show that for all $n\ge 0$ and 
any field $k$ the commutative algebras 
$k[x,y]/(x^2, y^{2+n})$ have representation dimension 3. 
For the proof, we give an explicit inductive construction
of a suitable generator-cogenerator. 
As a consequence, we obtain that the group algebras in characteristic 2 of 
the groups $C_2\times C_{2^m}$ have representation dimension 3.
Note that for $m\ge 3$ 
these group algebras have wild representation type.      
\medskip


\end{abstract}

\maketitle

\section{Introduction}

The representation dimension of an Artin algebra was defined by
M.\,Auslander \cite{Auslander} as a way of measuring homologically
how far an algebra is from being of finite representation type. 
(For the precise definition see Section \ref{Sec-defs} below.)
In fact, Auslander showed that the representation dimension of an
algebra is at most 2 if and only if the algebra has finite 
representation type, i.e. has only finitely many indecomposable
modules (up to isomorphism). 

For more than two decades it remained unclear whether Auslander's
philosophy and hopes were justified, one of the main problems 
being that the precise value of the representation dimension
of an algebra is very hard to determine. For instance, until
the late 1990's there were still only very few examples of algebras
known with representation dimension 3, and not a single algebra
with representation dimension $>3$. 

This situation changed dramatically over the last few years,
with spectacular recent progress which we briefly mention
here. O.\,Iyama proved that the representation dimension is always 
finite \cite{Iyama}. This integer attached to any algebra
is actually an invariant under various notions and levels of
equivalences of algebras. As it is defined, 
the representation dimension is invariant under Morita equivalences. 
In modern representation theory the notions of stable
equivalences and especially derived equivalences became 
increasingly important. It was recently shown by X.\,Guo \cite{Guo}
that the representation dimension is invariant under
stable equivalences. (For stable equivalences of Morita type this
was already shown by C.\,Xi \cite{Xi1}.) In particular, 
for selfinjective algebras the representation dimension is
invariant under derived equivalences, thus providing one of the
few explicit invariants of derived module categories known so far. 

A breakthrough was obtained by R.\,Rouquier \cite{RR-dims}
(but see also the new revised version \cite{RR-exterior}).
As a consequence of a deep general theory of dimensions of
triangulated categories he obtained the first examples of 
algebras with representation dimension $>3$, solving one
of the most fundamental open questions about this invariant. 
Much more than that, Rouquier obtained that the representation dimension
is unbounded, by showing that for any $n\ge 3$ 
the exterior algebra of an $n$-dimensional 
vector space has representation dimension $n+1$. 
Actually, the methods used there indicate a very geometric nature of 
the representation dimension in these cases, and from this approach
one can expect many more explicit examples in the near future. 

These recent results finally make it clear that Auslander's 
philosophy of measuring how far an algebra is from being 
of finite representation type is working; we get a new 
division of algebras according to the size of their
representation dimension. 
Note that this new division certainly does not run along
the lines of the classical division into finite, tame and wild representation
type. For instance, already Auslander showed that all finite-dimensional
hereditary algebras have representation dimension at most 3, not
depending on their representation type. 
There is something completely new to be discovered! So far, 
it remains mysterious in general what structural properties of an algebra
actually determine the representation dimension. 

In this paper, we will provide more examples of algebras 
of wild representation type having representation dimension 3. 
Our main motivation comes from group algebras. It is well known
that the representation type of a modular group algebra $kG$ of
a finite group over a field $k$ of characteristic $p$ 
is determined by a Sylow $p$-subgroup $P$ of $G$ 
\cite{BD}: $kG$ has finite
representation type if and only if $P$ is cyclic; $kG$ has tame
representation type precisely when $p=2$ and $P$ is dihedral,
semidihedral or quaternion. In all other cases $kG$ has wild 
representation type. Hence, in the cyclic case, 
the representation dimension of $kG$ is at most 2. 
It is also known that all group algebras of tame representation type
(and related algebras) 
have representation dimension 3. (Unfortunately, there is no
complete proof in the literature, but see \cite{Holm-AlgColl},
and \cite{Iyama2},4.4.2.)

From Rouquier's work one gets that for any $n\ge 3$ 
the mod 2 group algebra of an elementary abelian group 
$C_2\times \ldots\times C_2$ of order $2^n$ has representation dimension
$n+1$. One could hope that the geometric
methods used could be adapted and 
generalized to deal with elementary abelian $p$-groups in odd 
characteristic as well. However, these geometric methods do not work 
for elementary abelian groups of rank 2. 

The main aim of this paper is to give some partial results on the
rank 2 case, thus providing the first examples to fill this 
'rank-2-gap'. Note that the mod $p$ group algebra of 
$C_{p^m}\times C_{p^n}$ is isomorphic to the commutative ring
$k[x,y]/(x^{p^m}, y^{p^n})$. So the natural class of algebras to study
are the commutative rings $A_{n,m}:= k[x,y]/(x^m, y^n)$, for 
natural numbers $n,m\ge 2$. 
\smallskip

Our main result in this paper gives the
precise value of the representation dimension for the algebras
$A_{2,m}$. Note that this result is independent of the ground field. 

\begin{Theorem} \label{mainthm}
Let $k$ be a field and let $\Lambda _n=k[x,y]/(x^2, y^{2+n})$ 
for every non-negative integer $n$. Then the representation dimension 
of $\Lambda _n$ is 3. 
\end{Theorem}

As a consequence we get the following 
partial answer on the 'rank-2-gap' for group algebras.

\begin{Corollary} \label{cor-groups}
Let $k$ be a field of characteristic 2. Then for any $n\ge 1$ the 
group algebra $k(C_2\times C_{2^n})$ has representation dimension 3. 
\end{Corollary}

The paper is organized as follows. In Section \ref{Sec-defs} 
we recall the necessary definitions and 
background on representation dimensions. In particular, we outline
the method for determining the precise value of the
representation dimension. The final
Section \ref{Sec-main} then contains the
proof of the main result.


\section{Computing representation dimensions} \label{Sec-defs}

The representation dimension has been defined by M.\,Auslander 
\cite{Auslander} for
arbitrary Artin algebras. However, for the purpose of this paper,
by an algebra we mean a finite-dimensional algebra over a fixed 
field $k$. (It will turn out that our results will not depend 
on the ground field.) 

A finitely generated $A$-module $M$ is called a 
{\em generator-cogenerator} for $A$ if all projective indecomposable 
modules and all injective indecomposable $A$-modules occur 
as direct summands of $M$. 

Auslander's basic idea was to study homological properties 
of endomorphism rings of modules (instead of studying modules
directly). For any algebra $A$ the {\em representation dimension}
is defined to be 
$$\repdim(A):=\inf\{\gd(\End_A(M))\,|\,M\,\,
\mbox{generator-cogenerator}\}.$$
Note that this invariant a priori 
has values in $\mathbb{N}_0\cup\{\infty\}$. However, it was shown
by O.\,Iyama \cite{Iyama} that $\repdim(A)$ is always finite.

Auslander's fundamental result states that $\repdim(A)\le 2$ 
if and only if $A$ has finite representation type. More precisely,
$\repdim(A)=0$ if and only if $A$ is semisimple, and $\repdim(A)\neq
1$ for all algebras $A$. 

The definition indicates how to determine upper bounds for 
$\repdim(A)$; in fact, if one can show for some generator-cogenerator
$M$ that $\gd(\End_A(M))\le m$ then $\repdim(A)\le m$. In particular,
if one can find a generator-cogenerator $M$ with $\gd(\End_A(M))\le 3$,
and if $A$ is not of finite representation type, then 
one has shown that $\repdim(A)=3$. Of course, the notoriously difficult
problem is to find a suitable module $M$. 

However, this is exactly the strategy we will pursue 
successfully in the proof
of our main theorem. We will explain in Section
\ref{Sec-main} below how to construct inductively 
generator-cogenerators for the algebras $k[x,y]/(x^2,y^{2+n})$
with endomorphism rings of global dimension 3. 

We now describe the general method used later for
how to determine the global dimension
of the endomorphism ring of a generator-cogenerator $M$.

Let $M$ be a finitely generated module over a finite-dimensional 
algebra $A$.  The identity of $\ed_A(M)$ is the sum of the
``identity maps" on the indecomposable direct summands of $M$. 
Hence we have primitive idempotents of $\ed_A(M)$
corresponding to the summands of $M$. For any indecomposable 
summand $T$ of $M$ we denote the corresponding simple $\ed_A(M)$-module 
by $E_T$. The corresponding indecomposable projective 
$\ed_A(M)$-module $Q_T$ is given by all homomorphisms from $M$ 
to $T$, for abbreviation denoted by $Q_T=\Hm _A (M, T)=: (M, T)$. 
To prove that $\gd(\ed_A(M))\leq 3$, we explicitly construct a 
projective resolution 
with length $\leq 3$ for every simple module $E_T$. The general 
method is as follows. Recall that for a finitely generated 
$A$-module $M$, $\add M$ denotes the full subcategory
of $\mod A$ consisting of direct sums of direct summands of $M$.  
For any indecomposable summand $T$ of $M$, we first 
construct a suitable exact sequence $0\rightarrow K\rightarrow 
N_1 \rightarrow T$ with $N_1\in \add M$ with the following property:

\medskip
 (*) Every homomorphism from an indecomposable summand of $M$ to $T$, 
except the multiples of the identity on $T$, factors through $N_1$.

\medskip
 Applying the functor $(M, -)$, we get another short exact sequence 
$$0\rightarrow (M, K)\rightarrow (M, N_1) \rightarrow (M, T).$$ 
If the cokernel of $(M, N_1) \rightarrow (M, T)$ is 1-dimensional  
then the cokernel 
is $E_T$ (because $(M, T)$ is projective), i.e. we get the initial
part of a projective resolution of the simple $\End_A(M)$-module
$E_T$. 

If $K\in \add M$, then also $(M,K)$ is projective and 
we have constructed a projective 
resolution of $E_T$ of length 2. Otherwise, we construct another 
suitable short exact sequence
$0\rightarrow K'\rightarrow N_2 \rightarrow K$ with $N_2\in\add M$ 
with the following property:

\medskip
(**) Every map from an indecomposable summand of $M$ to $K$ factors 
through $N_2$. 

\medskip
 Applying the functor $(M, -)$ to this short exact sequence, we 
get a short exact sequence
$$0\rightarrow (M, K') \rightarrow (M, N_2) \rightarrow (M, K) 
\rightarrow 0.$$ 
If it happens that  $K'\in\add M$, we have $(M, N_1)$, $(M, N_2)$, 
$(M, K')$ and $(M, T)$ are all projective $\ed_A(M)$-modules. 
Hence we get a projective resolution
of the simple module $E_T$:
 $$0\rightarrow Q_{K'}\rightarrow Q_{N_2}\rightarrow Q_{N_1}
\rightarrow Q_T \rightarrow E_T\rightarrow 0$$
and $\pd E_T\leq 3$. 

The crucial aspect in this strategy is to come up with 
a generator-cogenerator $M$ for which this process really stops
at this stage,
i.e. for which $K'\in\add M$ in the second step above.


\section{Proof of the main theorem} \label{Sec-main}

Before proving our main theorem, we introduce some notations. 
For any non-negative integer $n$,
we use $\Lambda _n$ to denote the (commutative) 
algebra $k[x,y]/(x^2,y^{n+2})$. 
The quotient algebra 
$\Lambda _n/\soc (\Lambda _n)$ will be denoted by $A_n$.

Actually, in order to prove Theorem \ref{mainthm}
it suffices to show that $\repdim(A_n)=3$, due to the 
following general result.   
Recall that an algebra $A$ is called {\em basic} if all simple 
$A$-modules are of dimension 1. Note that our algebras
$k[x,y]/(x^2,y^{2+n})$ are basic since the trivial module $k$ is the 
only simple module.

\begin{Lemma}[\cite{EHIS}, Proposition 1.2]  \label{EHIS}
Let $\Lambda$ be a basic algebra, and let $P$ be an indecomposable 
projective-injective
$\Lambda$-module. 
Define $A=\Lambda/\soc (P)$. If $\rp(A)\leq 3$, then 
$\rp(\Lambda)\leq 3$. 
\end{Lemma}

The structure of the projective-indecomposable $A_n$-module
can be conveniently described diagrammatically as follows. 

\bigskip
\begin{center}
\begin{picture}(42,26)
\put(0,16){\includegraphics{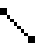}}
\put(8,8){\includegraphics{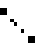}}
\put(16,0){\includegraphics{y}}
\put(0,24){\includegraphics{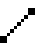}}
\put(8,24){\includegraphics{y}}
\put(8,16){\includegraphics{x}}
\put(16,16){\includegraphics{dottedy}}
\put(16,8){\includegraphics{x}}
\put(24,8){\includegraphics{y}}
\put(24,0){\includegraphics{x}}
\put(32,0){\includegraphics{y}}
{\footnotesize
\put(11,33){\includegraphics{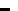}}
\put(14,29){\includegraphics{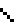}}
\put(18,25){\includegraphics{bracX}}
\put(22,24){\includegraphics{bracH}}
\put(24,22){\includegraphics{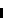}}
\put(27,24){$n$ sections}
\put(25,18){\includegraphics{bracX}}
\put(29,14){\includegraphics{bracX}}
\put(33,11){\includegraphics{bracV}}}
\end{picture}
\end{center}

Here, the vertices correspond to basis vectors of the module
and the edges describe the action of the generators $x$ and $y$
of $A_n$. 
\smallskip

Now, we use a table to introduce and illustrate the notations for 
several classes of $A_n$-modules which will be used 
throughout the paper. 

\medskip
\begin{center}
\begin{tabular}{c|c}
Notation  & Shape of the $A_n$-module \\
\hline
$U_i$, $0\leq i\leq n+1$  & 
\begin{picture}(40,30)
\put(0,16){\includegraphics{y}}
\put(8,8){\includegraphics{dottedy}}
\put(16,0){\includegraphics{y}}
\put(3,25){\includegraphics{bracH}}
\put(6,21){\includegraphics{bracX}}
\put(10,17){\includegraphics{bracX}}
\put(14,16){\includegraphics{bracH}}
\put(16,14){\includegraphics{bracV}}
\put(19,16){$i$ sections}
\put(17,10){\includegraphics{bracX}}
\put(21,6){\includegraphics{bracX}}
\put(25,3){\includegraphics{bracV}}
\end{picture} \\
\hline
$X$ & 
\begin{picture}(10,14)
\put(0,2){\includegraphics{x}}
\end{picture}\\
\hline
$A_i^j$, $i,j\geq 0, i+j\leq n$ & 
\begin{picture}(66,54)

\put(0,40){\includegraphics{y}}
\put(8,32){\includegraphics{dottedy}}
\put(16,24){\includegraphics{y}}
\put(24,16){\includegraphics{y}}
\put(32,8){\includegraphics{dottedy}}
\put(40,0){\includegraphics{y}}
\put(24,24){\includegraphics{x}}
\put(32,24){\includegraphics{y}}
\put(32,16){\includegraphics{x}}
\put(40,16){\includegraphics{dottedy}}
\put(40,8){\includegraphics{x}}
\put(48,8){\includegraphics{y}}
\put(48,0){\includegraphics{x}}
\put(56,0){\includegraphics{y}}
\put(3,49){\includegraphics{bracH}}
\put(6,45){\includegraphics{bracX}}
\put(10,41){\includegraphics{bracX}}
\put(14,40){\includegraphics{bracH}}
\put(16,38){\includegraphics{bracV}}
\put(19,40){$j$ sections}
\put(17,34){\includegraphics{bracX}}
\put(21,30){\includegraphics{bracX}}
\put(25,27){\includegraphics{bracV}}

\put(35,33){\includegraphics{bracH}}
\put(38,29){\includegraphics{bracX}}
\put(42,25){\includegraphics{bracX}}
\put(46,24){\includegraphics{bracH}}
\put(48,22){\includegraphics{bracV}}
\put(51,24){$i$ sections}
\put(49,18){\includegraphics{bracX}}
\put(53,14){\includegraphics{bracX}}
\put(57,11){\includegraphics{bracV}}

\end{picture}\\
\hline
$DA_i^j$, $i,j\geq 0, i+j\leq n$ &
\begin{picture}(66,62)

\put(0,48){\includegraphics{y}}
\put(8,40){\includegraphics{dottedy}}
\put(16,32){\includegraphics{y}}
\put(24,24){\includegraphics{y}}
\put(32,16){\includegraphics{y}}
\put(40,8){\includegraphics{dottedy}}
\put(48,0){\includegraphics{y}}
\put(32,24){\includegraphics{x}}
\put(40,24){\includegraphics{y}}
\put(40,16){\includegraphics{x}}
\put(48,16){\includegraphics{dottedy}}
\put(48,8){\includegraphics{x}}
\put(56,8){\includegraphics{y}}
\put(56,0){\includegraphics{x}}

\put(3,57){\includegraphics{bracH}}
\put(6,53){\includegraphics{bracX}}
\put(10,49){\includegraphics{bracX}}
\put(14,48){\includegraphics{bracH}}
\put(16,46){\includegraphics{bracV}}
\put(19,48){$j$ sections}
\put(17,42){\includegraphics{bracX}}
\put(21,38){\includegraphics{bracX}}
\put(25,35){\includegraphics{bracV}}

\put(43,33){\includegraphics{bracH}}
\put(46,29){\includegraphics{bracX}}
\put(50,25){\includegraphics{bracX}}
\put(54,24){\includegraphics{bracH}}
\put(56,22){\includegraphics{bracV}}
\put(59,24){$i$ sections}
\put(57,18){\includegraphics{bracX}}
\put(61,14){\includegraphics{bracX}}
\put(65,11){\includegraphics{bracV}}

\end{picture}\\
\hline
\end{tabular}
\end{center}
\bigskip

Note that the projective $A_n$-module will be denoted by $A_n^0$. 
By Lemma \ref{EHIS}, $\rp(A_n) \leq 3$ will 
imply that $\rp(\Lambda_n)\leq 3$. So we consider $A_n$ instead 
of $\Lambda _n$. 
For any $n\geq 0$, we will prove that  $\rp(A_n) \leq 3$ by 
constructing a generator-cogenerator $M_n$ for $A_n$
such that $\gd(\ed _{A_n}(M_n))\leq 3$. 

Now we give our constuction of the desired generator-cogenerators 
$M_n$ of $A_n$ as follows. 

For $n=0$ we put
$$M_0=A_0^0\oplus DA_0^0\oplus U_0\oplus U_1 \oplus X .$$
(Actually, for the algebra $A_0=k[x,y]/(x^2,y^2,xy)$ it is known
that $\repdim(A_0)=3$; in fact, $A_0$ is special biserial and
we can apply Corollary 1.3 from \cite{EHIS}. However, we give
the explicit generator-cogenerator to indicate the ``induction base''.) 

If $n$ is positive, we put
$$\displaystyle{M_n=\bigoplus_{i,j\geq 0, i+j\leq n} DA_i^j 
\oplus \bigoplus _{i\geq 0,j>0, i+j\leq n}}A_i^j 
\oplus A_n^0 \oplus \bigoplus_{0\leq i\leq n+1} U_i \oplus X$$

\noindent
\begin{Remark} {\em Let us point out more clearly why this is an
inductive construction, i.e. that one obtains $M_n$ from $M_{n-1}$
by a simple ``combinatorial'' procedure. Note that in $M_n$ we still 
have all summands of $M_{n-1}$, except the projective module $A_{n-1}^0$. 
The new summands in $M_n$ are $DA_i^{n-i}$ 
and $A_i^{n-i}$ (for $0\le i\le n$) and the uniserial module
$U_{n+1}$. With the exception of $A_n^0$ and its dual, they can 
be obtained by extending from summands of $M_{n-1}$ as follows. 
For all summands of $M_{n-1}$ we add an 
additional vertex on top (i.e. an additional basis vector of the module)
whenever possible; but we only allow this when no new squares are created
in the shape of the module. 
All modules obtained in this way are added as summands of $M_n$
(but only once, to avoid multiplicities). 

For instance, going from $M_0$ to $M_1$, the extension of
$DA_0^0$ to the four-dimensional module having the shape of a square 
is not allowed; on the other hand, the extension of $DA_0^0$ to the 
four-dimensional module $DA_0^1$ is allowed and gives a new summand. 
With this procedure, for $i=0,\ldots, n-1$ the new modules 
$DA_i^{n-i}$ are obtained from the summands $DA_i^{n-i-1}$ of $M_{n-1}$, 
the new summands $A_i^{n-i}$ are obtained from $A_i^{n-i-1}$
and $U_{n+1}$ is obtained from $U_n$. 
Note that the new projective module $A_n^0$ and the new injective 
module $DA_n^0$ can not be obtained with this method from summands
of $M_{n-1}$, so they have to be added separately. 

Hence, loosely speaking, the simple recipe for constructing $M_n$ from 
$M_{n-1}$ is as follows. 

(i) For all summands of $M_{n-1}$ construct all possible modules
obtained by adding a vertex on top, but without creating new squares
in the shapes of the modules. Add these modules as new summands
(if they are not already summands). 

(ii) Remove the old projective module $A_{n-1}^0$. (But keep all
other summands of $M_{n-1}$.)

(iii) Add the new projective $A_n^0$ and its dual $DA_n^0$ as summands. 
}
\end{Remark}

This inductive construction of a generator-cogenerator for $A_n$
is actually the crucial step in the proof of our main result. 
We will show below that the endomorphism ring of $M_n$ has
global dimension 3, which then implies that $A_n$ has representation
dimension 3, as claimed. 
\medskip

We are now in the position to give 
the proof our main result. 
\smallskip

\noindent  
{\bf Proof of Theorem \ref{mainthm}:} We have to show that $\gd 
(\ed _{A_n}(M_n))\leq 3$ for any $n\geq 0$.
\medskip

\noindent
{\bf (I) The case $n=0$.} 
We have $M_0=A_0^0\oplus DA_0^0\oplus U_0\oplus U_1 \oplus X$.  
\smallskip

     (1)The projective module $A_0^0$. Since the radical of the 
projective module $A_0^0$ is 
$U_0\oplus U_0 \in \add M_0$ and the exact sequence $0\rightarrow 
\rad A_0^0\rightarrow A_0^0$ clearly has the property (*), we can easily
get a projective resolution of the simple $\ed_{A_0}(M_0)$-module 
$E_{A_0^0}$:
$$0\rightarrow  Q_{U_0}\oplus Q_{U_0} \rightarrow Q_{A_0^0} 
\rightarrow E_{A_0^0}\rightarrow 0.$$
by applying the functor $(M_0, -)$ to $0\rightarrow \rad A_0^0
\rightarrow A_0^0$.  Hence $\pd E_{A_0^0}=1$. 
\smallskip

   (2) The injective module $DA_0^0$. There is an exact sequence 
$0\rightarrow U_0\rightarrow U_1\oplus X\rightarrow DA_0^0\rightarrow 0$ 
with the propertiy (*)   (Note
that there is no epimorphism from an indecomposable summand of $M_0$ 
to $DA_0^0$, except the multiples of the identity on $DA_0^0$). Since 
$U_0$ and 
the middle term are in $\add M_0$, by applying the functor $(M_0, -)$, we 
get a projective resolution of the simple $\ed_{A_0}(M_0)$-module 
$E_{DA_0^0}$:
   $$0\rightarrow Q_{U_0}\rightarrow Q_{U_1}\oplus Q_{X}\rightarrow 
Q_{DA_0^0}\rightarrow E_{DA_0^0}\rightarrow 0.$$
Hence $\pd E_{DA_0^0}\leq 2$. 
\smallskip

  (3) The module $X$. There is a short exact sequence $0\rightarrow 
U_0 \rightarrow A_0^0 \rightarrow X\rightarrow 0$ with the property (*) 
(Clearly, epimorphisms 
only occur from $A_0^0$ to $X$ and they factor through the middle term. 
Those non-epimorphisms from an indecomposable summand of $M_0$ to $X$ must
have image in $\rad X$ and clearly they also factor through the middle 
term.) Applying $(M_0, -)$ to the exact sequence, we get a 
projective resolution 
of the simple $\ed_{A_0}(M_0)$-module $E_{X}$:
  $$0\rightarrow Q_{U_0} \rightarrow Q_{A_0^0} \rightarrow Q_{X}
\rightarrow E_{X}\rightarrow 0$$
and hence $\pd E_{X}\leq 2$. 
\smallskip
 
  (4) The modules $U_i$ where $i=0,1$.  By symmetry, it suffices to show
that $\pd E_{U_0}\leq 2$. 
 For $U_0$, there is a short exact sequence  $0\rightarrow K \rightarrow 
DA_0^0 \oplus DA_0^0 \rightarrow U_0 \rightarrow 0$
with the property (*) and the kernel $K$ has the following shape

\smallskip
\begin{center}
\begin{picture}(34,10)
\put(0,0){\includegraphics{y}}
\put(8,0){\includegraphics{x}}
\put(16,0){\includegraphics{y}}
\put(24,0){\includegraphics{x}}
\end{picture}
\end{center}

Since $K$ is not a summand of $M_0$, we need to find a short exact 
sequence ending at $K$ with the property (**).  There is a short exact 
sequence 
$0\rightarrow U_0\oplus U_0 \rightarrow U_1\oplus A_0^0 \oplus X 
\rightarrow K \rightarrow 0$ with the property (**). Applying the 
functor $(M_0, -)$ to 
the above two short exact sequences and putting the resulting exact 
sequences together, we get a project resolution of $E_{U_0}$:
$$0\rightarrow Q_{U_0}\oplus Q_{U_0} \rightarrow Q_{U_1}\oplus 
Q_{A_0^0} \oplus Q_X\rightarrow Q_{DA_0^0} \oplus Q_{DA_0^0} 
\rightarrow Q_{U_0}
\rightarrow E_{U_0}\rightarrow 0$$
Hence $\pd E_{U_0}\leq 3$ . 
\smallskip

From the above, we have that $\gd (\ed M_0) \leq 3$, and 
as a consequence $\repdim(A_0)\leq 3$. Since $A_0$ is not
of finite representation type we can actually deduce that
$\repdim(A_0)= 3$. 
\medskip

\noindent
{\bf (II) The case $n>0$.} We have 
$$\displaystyle{M_n=\bigoplus _{i,j\geq 0, i+j\leq n} DA_i^j \oplus 
\bigoplus _{i\geq 0,j>0, i+j\leq n}}A_i^j 
\oplus A_n^0 \oplus \bigoplus_{0\leq i\leq n+1} U_i \oplus X.$$
Again, for any indecomposable summand $N$ of $M_n$,  we show that 
$\pd E_N\leq 3$ by explicitly constructing a projective resolution for the 
corresponding 
simple $\ed_{A_n}(M_n)$-module $E_N$. 
For the convenience of the reader, we list the indecomposable 
summands of $M_n$ as follows.

\begin{center}
\begin{tabular}{cccccc}
$A_n^0$ &  &  &  &  &  X   \\
  & $A_{n-1}^1$ &  $A_{n-2}^1$  & $\cdots$  & $A^1_1$    
  & $A_0^1$ \\
  & & $A_{n-2}^2$ & $\cdots$ & $A^2_1$ & $A_0^2$ \\
  & & & $\ddots$ & $\vdots$ & $\vdots$  \\
  & &  & & $A^{n-1}_1$  & $A_0^{n-1}$ \\
  & &  & &  & $A_0^n$   \\
\end{tabular}

\medskip

\begin{tabular}{ccccccc}
$DA_n^0$ & $DA_{n-1}^0$ &  $DA_{n-2}^0$ & $\cdots$   & $DA^0_1$      
 & $DA_0^0$     \\
 & $DA_{n-1}^1$ &  $DA_{n-2}^1$ & $\cdots$ & $DA^1_1$ 
 & $DA_0^1$     \\
 & &  $DA_{n-2}^2$ & $\cdots$ & $DA^2_1$ & $DA_0^2$     \\
 & & & $\ddots$  & $\vdots$ & $\vdots$ \\
 & & & & $DA^{n-1}_1$ & $DA_0^{n-1}$ \\
 & & & & & $DA_0^n$   \\
\end{tabular}

\medskip
$U_0, U_1, U_2, \cdots , U_n, U_{n+1}$

\end{center}

(1) The projective module $A_n^0$. Clearly, the exact sequence 
$0\rightarrow \rad A_n^0\rightarrow A_n^0\rightarrow 0$ has the 
property (*).  
Since $\rad A_n^0=A_{n-1}^1 \in \add M_n$,  applying the functor 
$(M_n, -)$ to the  exact sequence gives us a projective resolution of the 
simple $\ed_{A_n}(M_n)$ module $E_{A_n^0}$ of the form
$$0\rightarrow Q_{A_{n-1}^1}\rightarrow Q_{A_n^0} \rightarrow 
E_{A_n^0}\rightarrow 0. $$
Hence $\pd E_{A_n^0}=1$. 
\smallskip

(2) The module $X$.  There is a short exact sequence $0\rightarrow 
A_{n-1}^1  \rightarrow U_0\oplus A_n^0  \rightarrow X \rightarrow 0$ with
the property (*) (Actually, the only epimorphisms, apart from the 
multiples of the identity on $X$,  are from $A_n^0$, so clearly 
they factor through the middle term. All other maps factor through
the radical $\rad X=U_0$).  Since $A_{n-1}^1$ and the middle term 
are in $\add M_n$, applying the functor $(M_n, -)$, we get a 
projective resolution
of the simple $\ed_{A_n}(M_n)$-module $E_{X}$:
 $$0\rightarrow Q_{A_{n-1}^1}  \rightarrow Q_{U_0}\oplus Q_{A_n^0}  
\rightarrow Q_{X} \rightarrow E_{X}\rightarrow 0$$
Hence $\pd E_{X}\leq 2$. 
\smallskip

(3) The modules $A_i^1$.  If $i=0$, then there is a short exact 
sequence $0\rightarrow DA_0^1  \rightarrow  U_1\oplus U_0\oplus DA_1^0 
\rightarrow 
A_0^1\rightarrow 0$ with the property (*) (Except for the multiples of the 
identity on $A_0^1$, the epimorphisms are from $A_s^1$ ($s>0$) and 
$DA_t^0$. They 
all factor through $DA_1^0$).  Since all terms of the exact sequence 
are in $\add M_n$, by applying the functor $(M_n, -)$, we get a 
projective resolution 
of the simple $\ed_{A_n}(M_n)$-module $E_{A_0^1}$:
 $$0\rightarrow Q_{DA_0^1}  \rightarrow Q_{ U_1}\oplus Q_{ U_0}
\oplus Q_{ DA_1^0} \rightarrow Q_{ A_0^1}\rightarrow E_{A_0^1}\rightarrow 0$$
Hence $\pd E_{A_0^1}\leq 2$. 

If $i>0$, then there is a short 
exact sequence $0\rightarrow DA_i^1\rightarrow A_{i-2}^2\oplus 
DA_{i+1}^0\rightarrow A_i^1$
with the property (*) (Apart from the multiples of the identity 
on $A_i^1$, the epimorphisms to $A_i^1$ are from $A_s^1$ ($s>i$) 
and $DA_t^1$ ($t>i$). They
all factor through $DA_{i}^1$. The non-epimorphisms all factor 
through the middle term by the definition of the middle term).  
Note that all the terms of the 
short exact sequence are in $\add M_n$. Applying the functor 
$(M_n, -)$, we get a projective resolution of the simple 
$\ed_{A_n}(M_n)$-module $E_{A_i^1}$:
   $$0\rightarrow  Q_{DA_i^1} \rightarrow Q_{A_{i-1}^2}\oplus 
Q_{DA_{i+1}^0} \rightarrow Q_{A_i^1}  \rightarrow E_{A_i^1}\rightarrow 0$$
\smallskip

(4) The modules $A_i^j$ ($j>1, i>0$).  For every such module, 
there is a short exact sequence $0\rightarrow K \rightarrow 
A_{i-1}^{j+1}\oplus A_i^{j-1}\oplus 
DA_{i+1}^{j-1} \rightarrow A_i^j\rightarrow 0$ with the property (*) 
(In fact, the only epimorphisms from summands other than $A_i^j$ 
are from $A_s^j$ ($s>i$) 
and $DA_t^{j-1}$ ($t>i$). They all factor through $DA_{i+1}^{j-1}$. 
Other maps factor through $ A_{i-1}^{j+1}\oplus A_i^{j-1}$). 
The kernel $K$ has the following shape:

\begin{center}
{\footnotesize
 \begin{picture}(58,50)
\put(0,40){\includegraphics{y}}
\put(8,32){\includegraphics{dottedy}}
\put(16,24){\includegraphics{y}}
\put(24,16){\includegraphics{y}}
\put(32,8){\includegraphics{dottedy}}
\put(40,0){\includegraphics{y}}
\put(24,24){\includegraphics{x}}
\put(32,24){\includegraphics{y}}
\put(32,16){\includegraphics{x}}
\put(40,16){\includegraphics{dottedy}}
\put(40,8){\includegraphics{x}}
\put(48,8){\includegraphics{y}}
\put(48,0){\includegraphics{x}}
 \put(-11, 46){13}
\put(-8, 38){123} 
 \put(0, 30){123}
 \put(8, 22){123}
 \put(16, 14){123}
 \put(24, 6){123}
 \put(41, -2){3}
\put(60, 8){123}
\put(52, 16){123}
\put(44, 24){123}
\put(34, 32){23}
\put(3,49){\includegraphics{bracH}}
\put(6,45){\includegraphics{bracX}}
\put(10,41){\includegraphics{bracX}}
\put(14,40){\includegraphics{bracH}}
\put(16,38){\includegraphics{bracV}}
\put(19,45){$j$ sections}
\put(17,34){\includegraphics{bracX}}
\put(21,30){\includegraphics{bracX}}
\put(25,27){\includegraphics{bracV}}

\put(45,43){\includegraphics{bracH}}
\put(48,39){\includegraphics{bracX}}
\put(52,35){\includegraphics{bracX}}
\put(56,34){\includegraphics{bracH}}
\put(58,32){\includegraphics{bracV}}
\put(61,34){$i+1$ sections}
\put(59,28){\includegraphics{bracX}}
\put(63,24){\includegraphics{bracX}}
\put(67,21){\includegraphics{bracV}}

\end{picture}
}
\end{center}

Clearly, $K$ is not in $\add M_n$.  There is a short exact sequence 
$0\rightarrow DA_i^{j-1}\rightarrow A_{i-1}^j\oplus DA_i^j\oplus 
DA_{i+1}^{j-2}\rightarrow K\rightarrow 0$ with the property (**). 
Applying the functor $(M_n, -)$ to the above two exact sequences and putting
the resulting sequences together, we get a projective resolution 
of the simple $\ed_{A_n}(M_n)$-module $E_{A_i^j}$:

$$0\rightarrow Q_{DA_i^{j-1}} \rightarrow Q_{A_{i-1}^j}\oplus 
Q_{DA_i^j}\oplus Q_{DA_{i+1}^{j-2}}\rightarrow Q_{A_{i-1}^{j+1}}
\oplus Q_{A_i^{j-1}}\oplus 
Q_{DA_{i+1}^{j-1}}\rightarrow Q_{A_i^j}\rightarrow E_{A_i^j}\rightarrow 0$$
Hence $\pd E_{A_i^j}\leq 3$ for $j>1$ and $i>0$.
\smallskip

(5) The modules $A_0^j (j>1)$.  For any $A_0^j$ with $j>1$, there 
is a short exact sequence 
$0\rightarrow K  \rightarrow U_j\oplus A_0^{j-1}\oplus DA_1^{j-1} 
\rightarrow A_0^j\rightarrow 0$ with the property (*) (The only 
epimorphisms not from 
$A_0^j$ itself are from $A_s^j$ ($s>0$) or $DA_t^j$ ($t>0$). 
They all factor through $DA_1^{j-1}$). The kernel $K$ has the 
following shape:

\begin{center}
{\footnotesize
\begin{picture}(42,34)
\put(0,24){\includegraphics{y}}
\put(8,16){\includegraphics{dottedy}}
\put(16,8){\includegraphics{y}}
\put(24,0){\includegraphics{y}}
\put(24,8){\includegraphics{x}}
\put(32,8){\includegraphics{y}}
\put(32,0){\includegraphics{x}}
\put(-11,30){13}
\put(-8,22){123}
\put(0, 14){123}
\put(8,6){123}
\put(24, -2){3}
\put(35, 16){23}
\put(43,8){23}
\put(3,33){\includegraphics{bracH}}
\put(6,29){\includegraphics{bracX}}
\put(10,25){\includegraphics{bracX}}
\put(14,24){\includegraphics{bracH}}
\put(16,22){\includegraphics{bracV}}
\put(19,24){$j$ sections}
\put(17,18){\includegraphics{bracX}}
\put(21,14){\includegraphics{bracX}}
\put(25,11){\includegraphics{bracV}}
\end{picture}
}
\end{center}

Again, $K$ is not in $\add M_n$, but there is a short exact sequence 
$0\rightarrow U_j \rightarrow U_{j-1}\oplus U_{j+1} \oplus 
DA_1^{j-2} \rightarrow K 
\rightarrow 0$ with the property (**).  Applying the functor 
$(M_n, -)$ to the above two short exact sequences and putting 
together the resulting sequences,
we get a projective resolution of the simple $\ed_{A_n}(M_n)$-module 
$E_{A_0^j}$:
$$0\rightarrow Q_{U_j} \rightarrow Q_{U_{j-1}}\oplus Q_{U_{j+1}} 
\oplus Q_{DA_1^{j-2}}  \rightarrow Q_{U_j}\oplus Q_{A_0^{j-1}}
\oplus Q_{DA_1^{j-1}}
 \rightarrow Q_{A_0^j}\rightarrow E_{A_0^j}\rightarrow 0$$
  Hence $\pd E_{A_0^j}\leq 3$ for $j>1$. 
\smallskip

(6) The module $DA_n^0$.  There is a short exact sequence  
$0\rightarrow A_{n-1}^1\rightarrow DA_{n-1}^1\oplus A_n^0
\rightarrow DA_n^0\rightarrow 0$
with the property (*) and all terms of the short exact sequence 
are in $\add M_n$. Applying the functor $(M_n, -)$ to the short 
exact sequence, we get a 
projective resolution of the simple $\ed_{A_n}(M_n)$-module $E_{DA_n^0}$:
 $$0\rightarrow Q_{A_{n-1}^1}\rightarrow Q_{DA_{n-1}^1}\oplus 
Q_{A_n^0}\rightarrow Q_{DA_n^0}\rightarrow E_{DA_n^0}\rightarrow 0$$
Hence $\pd DA_n^0\leq 2$. 
\smallskip

(7) The modules $DA_i^0(0<i<n)$. For any $DA_i^0$ with $0<i<n$, 
there is a short exact sequence $0\rightarrow A_{i-1}^2
\rightarrow DA_{i-1}^1\oplus A_i^1
\rightarrow DA_i^0 \rightarrow 0$ with the property (*) 
(Except for the multiples of the identity on $DA_i^0$, the only 
epimorphisms to $DA_i^0$ are from 
$DA_s^0$ ($s>i$) and $A_t^1$ ($t\geq i$). They all factor 
through $A_i^1$).  Note that all terms of the short exact 
sequence are in $\add M_n$. Applying 
the functor $(M_n, -)$, we get a projective resolution of $E_{DA_i^0}$:
$$0\rightarrow Q_{A_{i-1}^2}\rightarrow Q_{DA_{i-1}^1}\oplus Q_{A_i^1} 
\rightarrow Q_{DA_i^0} \rightarrow E_{DA_i^0}\rightarrow 0$$
Hence we have $\pd E_{DA_i^0}\leq 2$ for $0<i<n$. 
\smallskip

(8) The module $DA_0^0$. There is a short exact sequence 
$0\rightarrow K\rightarrow X\oplus A_0^1\rightarrow DA_0^0
\rightarrow 0$ with the property (*) 
(Note that all epimorphisms from an indecomposable summand of 
$M_n$ to $DA_0^0$, except the multiples of the identity on 
$DA_0^0$, factor through 
$A_0^1$). The kernel $K$ is isomorphic to $A_0^0$, which is 
not in $\add M_n$. There is a short exact sequence $0\rightarrow 
A_{n-1}^1\rightarrow U_0\oplus U_0\oplus A_n^0\rightarrow K\rightarrow 0$ 
with the property (**). Applying the functor $(M_n, -)$ to the 
above two exact sequences and putting the resulting sequences 
together, we get a projective 
resolution of the simple $\ed M_n$ module $E-{DA_0^0}$:
$$0\rightarrow Q_{A_{n-1}^1}\rightarrow Q_{U_0}\oplus Q_{U_0}
\oplus Q_{A_n^0} \rightarrow Q_{X}\oplus Q_{A_0^1}\rightarrow 
Q_{DA_0^0}\rightarrow 
E_{DA_0^0}\rightarrow 0$$
Hence we have $\pd E_{DA_0^0}\leq 3$. 
\smallskip

(9) The modules $DA_i^j(i,j>0, i+j=n)$. For each module $DA_i^j$ 
with $i,j>0$ and $i+j=n$, there is a short exact sequence 
$0\rightarrow DA_{i-1}^j\rightarrow DA_{i-1}^{j+1}\oplus DA_i^{j-1} 
\rightarrow DA_i^j$ with the property (*)  and all the terms of the 
short exact sequence are 
in $\add M_n$. Applying the functor $(M_n, -)$ to the exact sequence, 
we get a projective resolution of the simple $\ed M_n$ module $E_{DA_i^j}$:
 $$0\rightarrow Q_{DA_{i-1}^j}\rightarrow Q_{DA_{i-1}^{j+1}}\oplus 
Q_{DA_i^{j-1}}\rightarrow Q_{DA_i^j}\rightarrow E_{DA_i^j}\rightarrow 0$$
Hence we have $\pd E_{DA_i^j}\leq 2$ for $i,j>0$ and $i+j=n$. 
\smallskip

(10) The modules $DA_0^j(j>0)$. First, we consider $DA_0^n$. There 
is a short exact sequence $0\rightarrow U_n\rightarrow U_{n+1}\oplus 
DA_0^{n-1}
\rightarrow DA_0^n$ with the property (*).  Applying the functor 
$(M_n, -)$ to the short exact sequence, we get a projective 
resolution of the simple 
$\ed_{A_n}(M_n)$-module $E_{DA_0^n}$:
$$0\rightarrow Q_{U_n}\rightarrow Q_{U_{n+1}}\oplus Q_{DA_0^{n-1}} 
\rightarrow Q_{DA_0^n}\rightarrow E_{DA_0^n}\rightarrow 0$$
Hence we have $\pd E_{DA_0^n}\leq 2$. For each $DA_0^j$ with $0<j<n$, 
there is a short exact sequence $0\rightarrow A_0^j\rightarrow DA_0^{j-1}
\oplus A_0^{j+1}\rightarrow DA_0^j\rightarrow 0$ with the property (*).  
Applying the functor $(M_n, -)$ to the short exact sequence, we get 
a projective 
resolution of the simple $\ed_{A_n}(M_n)$-module $E_{DA_0^j}$:
$$0\rightarrow Q_{A_0^j}\rightarrow Q_{DA_0^{j-1}}\oplus Q_{A_0^{j+1}}
\rightarrow Q_{DA_0^j}\rightarrow E_{DA_0^j}\rightarrow 0$$
Hence we have $\pd E_{DA_0^j}\leq 2$.
\smallskip

(11) The modules $DA_i^j (i,j>0, i+j<n)$. For each module $DA_i^j$ 
with $i,j>0$ and $i+j<n$, there is a short exact sequence 
$0\rightarrow K\rightarrow
DA_i^{j-1}\oplus DA_{i-1}^{j+1}\oplus A_i^{j+1}\rightarrow 
DA_i^j$ with the property (*) (Except the multiples of the identity 
on $DA_i^j$,  
the only epimorphisms to $DA_i^j$ are from $DA_s^j$ ($s>i$) and 
$A_t^{j+1}$ ($t\geq i$). They all factor through $A_i^{j+1}$. 
Other maps factor throgh the 
middle term by the definition of the middle term.) The kernel 
$K$ has the following shape:

\begin{center}
{\footnotesize
 \begin{picture}(58,50)
\put(0,40){\includegraphics{y}}
\put(8,32){\includegraphics{dottedy}}
\put(16,24){\includegraphics{y}}
\put(24,16){\includegraphics{y}}
\put(32,8){\includegraphics{dottedy}}
\put(40,0){\includegraphics{y}}
\put(24,24){\includegraphics{x}}
\put(32,24){\includegraphics{y}}
\put(32,16){\includegraphics{x}}
\put(40,16){\includegraphics{dottedy}}
\put(40,8){\includegraphics{x}}
\put(48,8){\includegraphics{y}}
\put(48,0){\includegraphics{x}}
\put(56,0){\includegraphics{y}}
\put(68, 0){3}
 \put(-11, 46){23}
\put(-8, 38){123} 
 \put(0, 30){123}
 \put(8, 22){123}
 \put(16, 14){123}
 \put(24, 6){123}
 \put(32, -2){123}
\put(60, 8){123}
\put(52, 16){123}
\put(44, 24){123}
\put(34, 32){13}

\put(3,49){\includegraphics{bracH}}
\put(6,45){\includegraphics{bracX}}
\put(10,41){\includegraphics{bracX}}
\put(14,40){\includegraphics{bracH}}
\put(16,38){\includegraphics{bracV}}
\put(19,45){$j+1$ sections}
\put(17,34){\includegraphics{bracX}}
\put(21,30){\includegraphics{bracX}}
\put(25,27){\includegraphics{bracV}}

\put(45,43){\includegraphics{bracH}}
\put(48,39){\includegraphics{bracX}}
\put(52,35){\includegraphics{bracX}}
\put(56,34){\includegraphics{bracH}}
\put(58,32){\includegraphics{bracV}}
\put(61,34){$i$ sections}
\put(59,28){\includegraphics{bracX}}
\put(63,24){\includegraphics{bracX}}
\put(67,21){\includegraphics{bracV}}
\end{picture}
}
\end{center}

Clearly, $K$ is not in $\add M_n$. However, there is a short exact 
sequence 
$0\rightarrow A_{i-1}^{j+1}\rightarrow DA_{i-1}^j \oplus A_i^j
\oplus A_{i-1}^{j+2}\rightarrow K\rightarrow 0$ with the property 
(**). Applying the 
functor $(M_n, -)$ to the above two exact sequences and putting 
the resulting sequences together, we get a projective resolution of the simple
$\ed_{A_n}(M_n)$-module $E_{DA_i^j}$:
$$0\rightarrow Q_{A_{i-1}^{j+1}}\rightarrow Q_{DA_{i-1}^j} 
\oplus Q_{A_i^j}\oplus Q_{A_{i-1}^{j+2}}\rightarrow 
Q_{DA_i^{j-1}}\oplus Q_{DA_{i-1}^{j+1}}\oplus Q_{A_i^{j+1}}
\rightarrow Q_{DA_i^j}\rightarrow E_{DA_i^j}\rightarrow 0$$
Hence we have $\pd E_{DA_i^j}\leq 3$ for $i,j>0$ and $i+j<n$. 
\smallskip

(12) The module $U_0$. There is a short exact sequence $0\rightarrow 
K\rightarrow DA_0^n\oplus DA_0\rightarrow U_0\rightarrow 0$ with the
property (*) (Except the multiples of the identity on $U_0$, all maps 
from an indecomposable summand of $M_n$ to $U_0$ factor either through
$DA_0^n$ or $DA_0$). The kernel $K$ has the following shape:

\begin{center}
{\footnotesize
\begin{picture}(42,26)
\put(0,16){\includegraphics{y}}
\put(8,8){\includegraphics{dottedy}}
\put(16,0){\includegraphics{y}}
\put(24,0){\includegraphics{x}}
\put(32,0){\includegraphics{y}}
\put(40,0){\includegraphics{x}}

\put(3,25){\includegraphics{bracH}}
\put(6,21){\includegraphics{bracX}}
\put(10,17){\includegraphics{bracX}}
\put(14,16){\includegraphics{bracH}}
\put(16,14){\includegraphics{bracV}}
\put(19,16){$n+1$ sections}
\put(17,10){\includegraphics{bracX}}
\put(21,6){\includegraphics{bracX}}
\put(25,3){\includegraphics{bracV}}
\end{picture}
}
\end{center}

Clearly, $K$ is not in $\add M_n$. However, there is a short exact 
sequence $0\rightarrow U_n\oplus U_0\rightarrow U_{n+1}\oplus A_0^n\oplus X
\rightarrow K\rightarrow 0$ with the property (**).  Applying the 
functor $(M_n,-)$ to the above exact sequence and putting the 
resulting sequences 
together, we get a projective resolution of the simple 
$\ed_{A_n}(M_n)$-module $E_{U_0}$:
$$0\rightarrow Q_{U_n}\oplus Q_{U_0}\rightarrow Q_{U_{n+1}}\oplus 
Q_{A_0^n}\oplus Q_{X} \rightarrow Q_{DA_0^n}\oplus Q_{DA_0}\rightarrow 
Q_{U_0}\rightarrow E_{U_0}\rightarrow 0$$
Hence we have $\pd E_{U_0}\leq 3$.
\smallskip

(13) The module $U_1$. There is a short exact sequence  $0\rightarrow 
K\rightarrow DA_0^1\oplus A_0^n\rightarrow U_1\rightarrow 0 $ with 
the property (*)
(Except the multiples of the identity on $U_1$, every epimorphism 
from an indecomposable summand of $M_n$ to $U_1$ factors through 
either $DA_0^1$ or
$A_0^n$. The maps having image in the radical of $U_1$ factor through 
$A_0^n$). The kernel $K$ has the following shape:

\begin{center}
{\footnotesize
\begin{picture}(50,34)
\put(0,24){\includegraphics{y}}
\put(8,16){\includegraphics{dottedy}}
\put(16,8){\includegraphics{y}}
\put(24,8){\includegraphics{x}}
\put(32,8){\includegraphics{y}}
\put(40,0){\includegraphics{y}}
\put(48,0){\includegraphics{x}}
\put(3,33){\includegraphics{bracH}}
\put(6,29){\includegraphics{bracX}}
\put(10,25){\includegraphics{bracX}}
\put(14,24){\includegraphics{bracH}}
\put(16,22){\includegraphics{bracV}}
\put(19,24){$n$ sections}
\put(17,18){\includegraphics{bracX}}
\put(21,14){\includegraphics{bracX}}
\put(25,11){\includegraphics{bracV}}
\end{picture}
}
\end{center}

There is another short exact sequence $0\rightarrow A_0^n\rightarrow 
U_n\oplus A_1^{n-1}\oplus DA_0^0\rightarrow K\rightarrow 0$ with 
the property (**). 
Applying the functor $(M_n, -)$ to the above two exact sequences and 
putting the resulting sequences together, 
we get a projective resolution of the simple 
$\ed_{A_n}(M_n)$-module $E_{U_1}$:

$$0\rightarrow Q_{A_0^n}\rightarrow Q_{U_n}\oplus Q_{A_1^{n-1}}\oplus 
Q_{DA_0^0}\rightarrow Q_{DA_0^1}\oplus Q_{A_0^n}\rightarrow Q_{U_1}
\rightarrow E_{U_1}\rightarrow 0$$
Hence we have $\pd E_{U_!}\leq 3$. 
\smallskip

(14) The modules $U_i, (1<i<n+1)$. For each module $U_i$ with $1<i<n+1$, 
there is a short exact sequence $0\rightarrow K\rightarrow U_{i-1}\oplus 
DA_0^i \oplus A_{i-1}^{n-i+1}\rightarrow U_i\rightarrow 0$ with the 
property (*) (In fact, the epimorphisms from an indecomposable 
summand of $M_n$ to 
$U_i$, except the multiples of the identity on $U_i$, either factor 
through $DA_0^{i}$ or $A_{i-1}^{n-i+1}$. Other maps factor through 
the $U_{i-1}$, which is 
the radical of $U_i$).  The kernel $K$ has the following shape:

\begin{center}
{\footnotesize
\begin{picture}(66,50)
\put(0,40){\includegraphics{y}}
\put(8,32){\includegraphics{dottedy}}
\put(16,24){\includegraphics{y}}
\put(24,16){\includegraphics{y}}
\put(32,8){\includegraphics{dottedy}}
\put(40,0){\includegraphics{y}}
\put(24,24){\includegraphics{x}}
\put(32,24){\includegraphics{y}}
\put(32,16){\includegraphics{x}}
\put(40,16){\includegraphics{dottedy}}
\put(40,8){\includegraphics{x}}
\put(48,8){\includegraphics{y}}
\put(48,0){\includegraphics{x}}
\put(56,0){\includegraphics{y}}
\put(64,-8){\includegraphics{y}}
\put(72,-8){\includegraphics{x}}

\put(3,49){\includegraphics{bracH}}
\put(6,45){\includegraphics{bracX}}
\put(10,41){\includegraphics{bracX}}
\put(14,40){\includegraphics{bracH}}
\put(16,38){\includegraphics{bracV}}
\put(19,43){$n-i+1$ sections}
\put(17,34){\includegraphics{bracX}}
\put(21,30){\includegraphics{bracX}}
\put(25,27){\includegraphics{bracV}}

\put(41,40){\includegraphics{bracH}}
\put(44,36){\includegraphics{bracX}}
\put(48,32){\includegraphics{bracX}}
\put(52,31){\includegraphics{bracH}}
\put(54,29){\includegraphics{bracV}}
\put(57,31){$i-1$ sections}
\put(55,25){\includegraphics{bracX}}
\put(59,21){\includegraphics{bracX}}
\put(63,18){\includegraphics{bracV}}

\put(-6, 45){3}
\put(2, 37){3}
\put(10, 29){3}
\put(18, 21){3}
\put(26, 13){3}
\put(34, 5){3}
\put(42, -3){3}

\put(35,31){23}
\put(43,23){123}
\put(51,15){123}
\put(59,7){123}
\put(67,-1){123}
\put(70,-15){2}
\put(84, -1){2}
\end{picture}
}
\end{center}

There is a short exact sequence $0\rightarrow A_{i-1}^{n-i+1}
\rightarrow A_i^{n-i}\oplus A_{i-2}^{n-i+2}\oplus DA_0^{i-1}
\rightarrow K\rightarrow 0$
with the property (**). Applying the functor $(M_n, -)$ to the 
above two exact sequences and putting together
the resulting sequences, 
we get a projective resolution of the simple 
$\ed M_n$ module $E_{U_i}$:
$$0\rightarrow Q_{A_{i-1}^{n-i+1}}\rightarrow Q_{A_i^{n-i}}
\oplus Q_{A_{i-2}^{n-i+2}}\oplus Q_{DA_0^{i-1}}\rightarrow Q_{U_{i-1}}\oplus 
Q_{DA_0^i} \oplus Q_{A_{i-1}^{n-i+1}}\rightarrow Q_{U_i}
\rightarrow E_{U_i}\rightarrow 0$$
Hence $\pd E_{U_i}\leq 3$ for $1<i<n+1$. 
\smallskip

(15)The module $U_{n+1}$. There is a short exact sequence 
$0\rightarrow A_{n-1}^1\rightarrow U_n\oplus A_n^0\rightarrow 
U_{n+1}\rightarrow 0$ with 
the property (**) (In fact, except the multiples of the identity 
on $U_{n+1}$, the only epimorphisms are from $A_n^0$. They factor 
through $A_n^0$. Other maps
factor through the radical $U_n$). Since all terms of the short 
exact sequence are in $\add M_n$, by applying the functor $(M_n, -)$, 
we get a projective 
resolution of the simple $\ed_{A_n}(M_n)$-module $E_{U_{n+1}}$:
$$0\rightarrow Q_{A_{n-1}^1}\rightarrow Q_{U_n}\oplus Q_{A_n^0}
\rightarrow Q_{U_{n+1}}\rightarrow E_{U_{n+1}}\rightarrow 0$$
Hence we have $\pd E_{U_{n+1}}\leq 2$. 
\smallskip

Thus we have shown that $\gd \ed _{A_n}(M_n)\leq 3$ for any 
non-negative integer $n$.  Hence we have $\rp A_n\leq 3$ for 
all $n\geq 0$. On the other hand, the algebras 
$A_n$ ($n\geq 0$) are not representation finite, so we have 
$\rp A_n>2$ for all $n\geq 0$. Altogether, we come to the 
conclusion that $\rp A_n=3$ for all non-negative
integers $n$ and this completes the proof of Theorem \ref{mainthm}.  
$\hfill\diamond$



\vskip1cm

\noindent
{\sc Thorsten Holm
\smallskip

Department of Pure Mathematics, University of Leeds

Leeds LS2 9JT, England}
\smallskip

E-mail:
{\tt tholm@maths.leeds.ac.uk}

\vskip.5cm

\noindent
{\sc Wei Hu
\smallskip

Department of Pure Mathematics, University of Leeds

Leeds LS2 9JT, England

and

School of Mathematical Sciences, Beijing Normal University

100875 Beijing, P. R. China}
\smallskip

E-mail:
{\tt huwei@maths.leeds.ac.uk}


\begin{thebibliography}{99}

\bibitem{Auslander}{{\sc M.\,Auslander}, 
Representation dimension of Artin algebras. Queen Mary College,
Mathematics Notes, University of London. 1971.
Also in: I.\,Reiten, S.\,Smal{\o}, {\O}.\,Solberg (Eds.),
Selected works of Maurice Auslander, Part I, Amer. Math. Soc.,
Providence, RI, 1999, pp. 505-574.}

\bibitem{BD}
{{\sc V.M. Bondarenko, J.A. Drozd},
The representation type of finite 
groups. (Russian) Zap. Nau\v cn. Sem. Leningrad. Otdel. Mat. Inst. 
Steklov. (LOMI)  71  (1977), 24-41. English translation:
J. Soviet Math. 20 (1982), 2515-2528.}

\bibitem{EHIS}{{\sc K.\,Erdmann, T.\,Holm, O.\,Iyama and 
J.\,Schr\"oer}, Radical embeddings and representation dimension. 
Advances in Math. 185 (2004), 159-177.}

\bibitem{Guo}{{\sc X.\,Guo}, 
Representation dimension: An invariant under stable equivalence. 
Trans. Amer. Math. Soc. 357 (2005) 3255-3263.}

\bibitem{Holm-AlgColl}{{\sc T.\,Holm}, Representation dimension of some 
tame blocks of finite groups. Algebra Colloq. 10 (2003), 275-284.}

\bibitem{Iyama}{{\sc O.\,Iyama}, Finiteness of representation dimension.
Proc.\,Amer.\,Math.\,Soc. 131\,(2003), 1011-1014.}
 
\bibitem{Iyama2}
{{\sc O.\,Iyama}, Rejective subcategories of artin algebras and orders.
Preprint (2003), math.RT/0311281.}

\bibitem{RR-dims} {{\sc R.\,Rouquier}, Dimensions of 
triangulated categories. Preprint (2003), math.CT/0310134}

\bibitem{RR-exterior}{{\sc R.\,Rouquier}, Representation dimension 
of exterior algebras. Preprint (2005), available at
{\tt http://www.math.jussieu.fr/$\sim$rouquier/
preprints/preprints.html}
}

\bibitem{Xi1}{{\sc C.\,Xi}, Representation dimension and 
quasi-hereditary algebras. Advances in Math. 168 (2002), 193-212.}

\end{thebibliography}
\end{document}